\documentclass[a4paper,11pt]{article}

\usepackage{microtype}
\usepackage[]{algorithm2e}
\usepackage{chngpage}

\usepackage[T1]{fontenc}
\usepackage[utf8]{inputenc}

\usepackage{amsfonts}
\usepackage{amsmath}
\numberwithin{equation}{section}

\usepackage{amssymb}
\usepackage{mathtools}
\usepackage[binary-units=true]{siunitx}

\usepackage{bm}
\usepackage{breqn}

\usepackage{amsthm}
\newtheorem{theorem}{Theorem} 
\newtheorem{lemma}[theorem]{Lemma}

\usepackage{caption}
\usepackage{subcaption}
\usepackage{diagbox}

\usepackage{txfonts}
\usepackage{sansmath}

\usepackage{geometry}
\newgeometry{vmargin={25mm}, hmargin={25mm,25mm}}

\usepackage{multicol}
\usepackage{booktabs}
\usepackage{graphicx}
\usepackage{float}

\makeatletter
\newcommand\tabfill[1]{%
\dimen@\linewidth%
\advance\dimen@\@totalleftmargin%
\advance\dimen@-\dimen\@curtab%
\parbox[t]\dimen@{#1\ifhmode\strut\fi}%
}

\usepackage{dcolumn}
\newcolumntype{d}[1]{D{.}{.}{#1}}

\usepackage{tikz}
\usetikzlibrary{shapes,arrows}
\usetikzlibrary{matrix}
\usepackage{pgfplots}
\usepackage{pgfplotstable}
								
\usepackage{color}
\definecolor{thered}{rgb}{0.65,0.04,0.07}
\definecolor{thegreen}{rgb}{0.06,0.44,0.08}
\definecolor{theblue}{rgb}{0.02,0.2,0.68}
\definecolor{clr1}{HTML}{bdd7e7}
\definecolor{clr2}{HTML}{6baed6}
\definecolor{clr3}{HTML}{3182bd}
\definecolor{clr4}{HTML}{08519c}

\renewcommand{\@algocf@capt@plain}{above}

\usepackage[unicode=true,
            bookmarks=true,
            bookmarksnumbered=false,
            bookmarksopen=false,
            breaklinks=true,
            pdfborder={0 0 1},
            backref=false,
            colorlinks=true,
            linkcolor=black,
            urlcolor=theblue,
            citecolor=theblue,
            hyperfootnotes=false]
           {hyperref}
\hypersetup{pdftitle={},
            pdfauthor={B. C. Vermeire}}

\begin{document}

\title{Energy Conservative Relaxation-Free Runge-Kutta Schemes}

\author{Mohammad R. Najafian and Brian C. Vermeire\\\\
  \textit{Department of Mechanical, Industrial, and Aerospace Engineering}\\
	\textit{Concordia University}\\
	\textit{Montreal, QC, Canada}}
\maketitle

\begin{abstract}
A wide range of physical phenomena exhibit auxiliary admissibility criteria, such as conservation of entropy or various energies, which arise implicitly under exact solution of their governing PDEs. However, standard temporal schemes, such as classical Runge-Kutta (RK) methods, do not enforce these constraints, leading to a loss of accuracy and stability. Previously, the Incremental Directional Technique RK (IDT-RK) and Relaxation Runge-Kutta (R-RK) approaches have been proposed to address this. However, these lead to a loss of accuracy in the case of IDT-RK, or a loss of step size control in the case of R-RK. In the current work we propose Relaxation-Free Runge-Kutta (RF-RK) schemes, which conserve energy, maintain order of accuracy, and maintain a constant step size, alleviating many of the limitations of the aforementioned techniques. Importantly, they do so with minimal additional computational cost compared to the base RK scheme. Numerical results demonstrate that these properties are observed in practice for a range of applications. Therefore, the proposed RF-RK framework is a promising approach for energy conservative time integration of systems of PDEs.

\end{abstract}

\section{Introduction}

Many physical systems yield measurable quantities that either remain constant or evolve monotonically over time. Examples of such behavior include explicitly conserved quantities, such as mass, momentum, and energy, or implicitly conserved quantities such as entropy in the compressible Euler equations with smooth solutions. When solving these Partial Differential Equations (PDEs) approximately, an important indicator of the quality of the numerical solution is preserving physical progression of these parameters \cite{hairerGeometricNumericalIntegration2006}. Often, these constraints are enforced explicitly through chosen conservation laws, such as conservation of mass, momentum, and total energy with the Euler equations.However, in some cases, such as entropy, they are enforced indirectly under exact evolution of conservation laws. Failure to maintain these at approximate discrete solutions can result in non-physical solutions \cite{GEAR1992303}, stability issues for long-time integration \cite{arakawaComputationalDesignLongTerm1997}, or increased numerical error
 \cite{defrutosAccuracyConservationProperties1997, duranNumericalIntegrationRelative2000,ranochaRateErrorGrowth2021} .Therefore, numerical schemes that can provide structure-preserving solutions, maintaining these implicitly conserved quantities, are of great importance.
  
Solving time-dependent PDEs commonly involves two main steps: first, the domain undergoes spatial discretization, transforming the system into semi-discrete time-dependent Ordinary Differential Equations (ODEs). Subsequently, these ODEs are solved using a temporal integration scheme to obtain a numerical solution that is fully-discrete  in space and time. Each of these  steps can potentially contaminate conserved quantities. For spatial discretization, there is a rich literature on stability-preserving techniques concerning kinetic energy and entropy for the Euler and Navier-Stokes equations, see e.g. \cite{chenReviewEntropyStable2020, chanDiscretelyEntropyConservative2018a, delreyfernandezReviewSummationbypartsOperators2014, kuyaHighorderAccurateKineticenergy2021}. Nevertheless, the resulting semi-discrete ODEs should be coupled with a time integration scheme to proceed in time, and the important question remains: if the conserved quantities are maintained after numerical integration. To answer that, stand-alone integration schemes can be studied in terms of conserving nonlinear stability properties of the ODE system. Here we consider the case of energy conservation.

To proceed, consider the following ODE system, which can be a semi-discrete representation of an initial PDE problem discretized via a spatial discretization

\begin{subequations}\label{IVP}
  \begin{align}
    u'(t) &= f(t,u(t)), \\
    u(t_0) &= u_0,
  \end{align}
\end{subequations} 
where $u(t)$ is a $m\times 1$ real vector of a solution in Hilbert space. A common quantity of interest is the $L_2$ inner product norm of the solution vector, or herein called energy: $\langle u,u \rangle= ||u||^2$. For the ODE system \ref{IVP}, energy is a continuous function of time: $||u(t)||^2$, while for its numerical solution, energy is only available at each discrete time step $n$: $||u(t_0+ n\Delta t)||^2= ||u^n||^2$. We call the semi-discrete system \ref{IVP} energy \textit{dissipative} (or strongly stable) if it satisfies an energy decay relation \cite{ketchesonRelaxationRungeKutta2019, sunStrongStabilityExplicit2019}

\begin{equation}\label{dissip_prob}
\frac{d}{dt}||u(t)||^2=2 \langle u, f(t,u) \rangle \leq 0.
\end{equation}
Maintaining this decay in energy by an integration scheme implies that the energy at each time step should not be increasing

\begin{equation}\label{dissip_sol}
||u^{n+1}||^2 \leq ||u^n||^2.
\end{equation} \\
Time integration schemes that guarantee non-increasing energy at each time step for dissipative problems are called \textit{monotonicity preserving} (also known as strong stability preserving) integration methods.

Moreover, system \ref{IVP} will be called energy \textit{conservative} if the time rate change of energy is exactly zero
 
\begin{equation}\label{conserv_prob}
\frac{d}{dt}||u(t)||^2=2 \langle u, f(t,u) \rangle = 0.
\end{equation}
Such integration schemes that, for conservative systems, maintain  energy unchanged at each step up to machine precision, are called energy $conservative$ integration schemes \cite{ketchesonRelaxationRungeKutta2019}, yielding

\begin{equation}\label{conserv_sol}
||u^{n+1}||^2 = ||u^n||^2.
\end{equation} \\

As previously mentioned, conservation and monotonicity preservation play a pivotal role in securing accurate and physically meaningful solutions \cite{biswasMultipleRelaxationRungeKutta2023, sunStabilityFourthOrder2017}. However, many widely used integration schemes fall short of guaranteeing these properties. For dissipative systems, it has been shown that many explicit schemes cannot preserve monotonicity, even for linear autonomous systems. For instance, Sun \& Shu \cite{sunStabilityFourthOrder2017} demonstrated that the classical fourth order RK method may adversely result in an increase in energy at the first time step, no matter how small the step size is. Moreover, no classical explicit RK method can preserve energy up to machine precision for conservative systems. They all lead to spurious energy change above their truncation error in the general case \cite{lozanoEntropyProductionExplicit2018}.

There have been attempts to modify explicit RK schemes to make them energy conservative while maintaining favorable accuracy and efficiency properties. In this regard, Calvo and coauthors \cite{calvoPreservationInvariantsExplicit2006} developed a directional projection technique, in which at the end of each time step, the solution is projected along an oblique direction to the conservative manifold. This technique has been further developed by Ketcheson \cite{ketchesonRelaxationRungeKutta2019} and referred to as the Relaxation Runge-Kutta technique (R-RK). With this method, the resulting explicit scheme is energy conservative, and the order of accuracy is retained. However, it comes at the expense of step size relaxation at each step. If instead one desires to keep the step size unchanged, it can be done with the so-called Incremental Directional Technique (IDT-RK), but the order of accuracy of the scheme will be decreased by one \cite{calvoPreservationInvariantsExplicit2006}. 

To address the limitations of both R-RK and IDT-RK techniques, this paper introduces a novel approach making explicit RK methods energy conservative, while the step size remains constant and the order of accuracy is retained. We refer to this method as the \textit{Relaxation-Free} Runge-Kutta approach (RF-RK). This method's additional computational cost at each step remains minimal, and is comparable to that of R-RK and IDT-RK methods. 

The rest of this paper is organized as follows.
In Section \ref{Sec_Methodology} a brief description of explicit Runge-Kutta schemes is provided. It then reviews relaxation Runge-Kutta methods, and the relaxation-free method is introduced. In Section \ref{Sec_Acc_Stab} it will be demonstrated that the accuracy of the original RK scheme is preserved with the proposed relaxation-free technique. Then, the linear stability regions of RF-RK schemes compared to their base RK method will be studied. Finally, numerical examples and comparisons between these energy conserving methods and unmodified RK schemes have been provided in Section \ref{Sec_exs}. Finally, conclusions and recommendations for future work are presented in Section \ref{Sec_conclusion}.

\section{Methodology} \label{Sec_Methodology}

This section describes classical RK methods and the two energy conserving techniques: the previous R-RK approach, and the proposed RF-RK approach. Section \ref{Subs_Classical_RK} briefly describes the classical RK framework, along with the energy change introduced at each step. Section \ref{Subs_RRK} briefly reviews how the relaxation technique is defined to cancel out the spurious energy change. Then, in Section \ref{Subs_RFRK} , our alternative technique, called a relaxation-free method, is presented. 

\subsection{Classical Runge-Kutta Schemes} \label{Subs_Classical_RK}

Given a solution vector at the time step $n$, denoted as $u^n=u(t_n)$, the next-step solution using a classical $s$ stage Runge-Kutta scheme is approximated as

\begin{subequations}\label{RK}
  \begin{align}
    y_i &= u^n + \Delta t \sum_{j=1}^s a_{ij} f_j , \quad i=1, ...,s,\label{RK_1} \\
    u(t_n+\Delta t) \approx u^{n+1} &= u^n + \Delta t \sum_{j=1}^s b_j f_j. \label{RK_2}
  \end{align}
\end{subequations} 
where $f_j$ is the \textit{j}th stage derivative
 
\[  f_j = f(t_n+ c_j\Delta t, y_j), \]
and the coefficients $a_{ij}$, $b_j$, and $c_j= \sum_i a_{ij}$ are elements of the scheme's Butcher tableau $A$, $b$, and $c$, respectively. An RK scheme is explicit when its matrix $A$ is lower triangular.

Following Ketcheson \cite{ketchesonRelaxationRungeKutta2019}, the change in energy after each time step will be

\begin{align*}
||u^{n+1}||^2 - ||u^n||^2 &= \bigg\Vert u^n + \Delta t \sum_{j=1}^s b_j f_j \bigg\Vert ^2 - ||u^n||^2, \\
&= 2\Delta t \sum_{j=1}^s b_j \langle u^n, f_j \rangle + \Delta t^2 \sum_{i,j=1}^s b_i b_j \langle f_i, f_j \rangle, \\
&= 2\Delta t \sum_{j=1}^s b_j \langle y_j, f_j \rangle + 2\Delta t \sum_{j=1}^s b_j \langle u^n - y_j, f_j \rangle + \Delta t^2 \sum_{i,j=1}^s b_i b_j \langle f_i, f_j \rangle.
\end{align*}
By using Eq. (\ref{RK_1}), this can be written as

\begin{equation} \label{RK_energy_change}
||u^{n+1}||^2 - ||u^n||^2 = 2\Delta t \sum_{j=1}^s b_j \langle y_j, f_j \rangle - 2\Delta t^2 \sum_{i,j=1}^s b_ia_{ij} \langle f_i, f_j \rangle + \Delta t^2 \sum_{i,j=1}^s b_ib_j \langle f_i, f_j \rangle .
\end{equation}
The first term on the right-hand side of Eq. (\ref{RK_energy_change}) depends on the spatial semi-discretization, and will be zero for conservative systems, or negative for dissipative ones. On the other hand, the summation of the remaining two terms, denoted by $\Delta E_t^n$

\begin{equation}
\Delta E_t^n= - 2\Delta t^2 \sum_{i,j=1}^s b_ia_{ij} \langle f_i, f_j \rangle + \Delta t^2 \sum_{i,j=1}^s b_ib_j \langle f_i, f_j \rangle,
\end{equation}
is the spurious energy introduced by numerical integration. If this is zero, then energy conservation, or monotonicity, of the semi-discrete system will be preserved. However, when $\Delta E_t^n$ deviates from zero, which is possible for all explicit RK schemes, energy conservation can be lost, and monotonicity may not be preserved \cite{ketchesonRelaxationRungeKutta2019}. The objective is to enforce $\Delta E^n_t =0$ at each step by modifying either the parameters of the integration scheme or the time step size at each step in Eq. (\ref{RK_2}).

\subsection{Relaxation Runge-Kutta Schemes} \label{Subs_RRK}

The relaxation family of Runge-Kutta schemes of Ketcheson \cite{ketchesonRelaxationRungeKutta2019} is created by replacing Eq. (\ref{RK_2}) with the following expression

\begin{equation}\label{RRK}
u(t_n+ \gamma_n\Delta t) \approx u^{n+1} = u^n + \Delta t \gamma_n\sum_{j=1}^s b_j f(t_n+ c_j\Delta t, y_j).
\end{equation}
This means that while intermediate solutions and their derivatives in Eq. (\ref{RK_1}) are computed using the original step size, $\Delta t$, the next step solution in Eq. (\ref{RRK}) in fact falls at $t_n + \gamma_n \Delta t$. This inherent inconsistency in the time step size provides freedom to have control over energy change at each step by fine-tuning the so-called relaxation variable $\gamma_n$.

With this relaxation method, the energy change at each step becomes

\begin{equation}\label{energy_RRK}
\begin{split}
||u_{\gamma}^{n+1}||^2 - ||u^n||^2 = 2\gamma_n\Delta t \sum_{j=1}^s b_j \langle y_j, f_j \rangle &- 2\gamma_n\Delta t^2 \sum_{i,j=1}^s b_ia_{ij} \langle f_i, f_j \rangle \\
&+ \gamma_n^2\Delta t^2 \sum_{i,j=1}^s b_ib_j \langle f_i, f_j \rangle.
\end{split}
\end{equation}
The last two terms, which comprise $\Delta E_t^n$, can be eliminated by solving for $\gamma_n$ via

\begin{equation}\label{gamma_n}
\gamma_n= \left\{ \begin{array}{lcl}
1  & ||\sum_{j=1}^s b_jf_j||^2=0 \\
 \frac{2 \sum_{i,j=1}^s b_ia_{ij} \langle f_i, f_j \rangle}{\sum_{i,j=1}^s b_ib_j\langle f_i, f_j \rangle } & ||\sum_{j=1}^s b_jf_j||^2 \neq 0 \end{array} \right.
\end{equation}

Ketcheson \cite{ketchesonRelaxationRungeKutta2019} demonstrated that while being energy conservative, this method retains the order of accuracy of the original RK scheme, and the additional cost is limited to the inexpensive calculation of $s+1$ inner products. However, the drawback of this method persists: by relaxing the time step size, the user will not have direct control over the actual step sizes taken. Opting for a fixed step size, denoted as the IDT-RK approach, means that the resulting integration scheme has parameters $\gamma_n b_j$ instead of $b_j$. This creates a local inconsistency since in general, $\sum_j^s \gamma_n b_j \neq 1$, which leads to reduced accuracy. Therefore, the user has to decide between relaxing time steps or losing an order of accuracy.

The following section introduces our proposed technique, called here a relaxation-free RK scheme (RF-RK), enabling explicit RK schemes	 to be energy conservative without step size relaxation or decreasing the scheme's order of accuracy.

\subsection{Relaxation-Free Runge-Kutta Schemes} \label{Subs_RFRK}

We introduce a family of energy conservative RK schemes, called here relaxation-free RK methods, as they do not impose relaxation on the time step size. With this method, the coefficients $b_j$ of the original RK scheme in Eq. (\ref{RK_2}) are substituted with modified values, $\hat{b}_j$

\begin{equation}\label{RFRK}
u(t_n+\Delta t) \approx u^{n+1} = u^n + \Delta t \sum_{j=1}^s \hat{b}_j f_j,
\end{equation}

\begin{equation}\label{hat_b}
\hat{b}_j= b_j + k_j \epsilon_n,
\end{equation}
where $\epsilon_n$ is a real-valued parameter to be calculated at each step $n$ to enforce elimination of spurious energy production or dissipation induced by numerical integration, and $k_j$ are constant multipliers that must adhere to certain conditions to follow. First, similar to all consistent RK schemes that satisfy $\sum_j b_j=1$, this requirement needs to be followed by the new $\hat{b}_j$ parameters: $\sum_j \hat{b}_j=1$. With Eq. (\ref{hat_b}), this can be guaranteed if $\sum_j k_j=0$. Moreover, it will be shown in Section \ref{Subs_Acc} that to keep the order of accuracy unchanged, these $k_j$ parameters need to be selected such that $\sum_j k_jc_j \neq 0$. So, to design a RF-RK method from a base RK scheme, we need the following

\begin{subequations}\label{k_conditions}
\begin{align}
\sum_{i=1}^s k_i &=0, \label{sum_ki} \\
\sum_{i=1}^s k_i c_i &\neq 0 . \label{sum_kici}
\end{align}
\end{subequations}

With this modified scheme, energy variation at each step takes the form:

\begin{equation} \label{RFRK_energy_change}
\begin{split}
||u^{n+1}||^2 - ||u^n||^2 = 2\Delta t \sum_{j=1}^s (b_j+k_j\epsilon_n) \langle y_j, f_j \rangle &- 2\Delta t^2 \sum_{i,j=1}^s (b_i+k_i\epsilon_n)a_{ij} \langle f_i, f_j \rangle \\ &+ \Delta t^2 \sum_{i,j=1}^s (b_i+k_i\epsilon_n)(b_j+k_j\epsilon_n) \langle f_i, f_j \rangle.
\end{split}
\end{equation}
Again, the summation of the last two terms is the numerically-induced energy change. Aiming to make this zero, a quadratic equation has to be solved to find the appropriate value for $\epsilon_n$

\begin{equation}\label{epsilon_eq}
A_n^* \epsilon_n^2 + B_n^*\epsilon_n + C_n^*=0,
\end{equation}
where $A^*$, $B^*$, and $C^*$ are

\begin{subequations}\label{ABC}
\begin{align}
A_n^* &= \sum_{i,j=1}^s k_i k_j \langle f_i, f_j \rangle, \label{A} \\
B_n^* &= -2 \sum_{i,j=1}^s k_i a_{ij} \langle f_i, f_j \rangle + 2 \sum_{i,j=1}^s k_i b_j \langle f_i, f_j \rangle, \label{B} \\
C_n^* &= -2 \sum_{i,j=1}^s b_i a_{ij} \langle f_i, f_j \rangle +  \sum_{i,j=1}^s b_i b_j \langle f_i, f_j \rangle. \label{C}
\end{align}
\end{subequations}
For a special case of constant stage derivatives, $f_i=f$, which happens in steady state solutions, we would have $A_n^*=(\sum_{i} k_i)^2 \langle f,f \rangle =0$, and $C_n^*= -2 \sum_{i} b_ic_i \langle f,f \rangle  + (\sum_{i} b_i)^2 \langle f,f \rangle =0$ for all RK schemes of second order or higher. In this case, the solution for Eq. (\ref{epsilon_eq}) becomes $\epsilon_n=0$, which aligns intuitively with the notion that in a steady state solution, there is no change in energy with the original RK scheme. However, when the condition $f_i =f$ does not hold, the solution depends on the sign of parameter $\hat{\Delta}_n$

\begin{equation}\label{Delta}
\hat{\Delta}_n= (B_n^*)^2 - 4A_n^*C_n^*.
\end{equation}
When $\hat{\Delta}_n$ takes a positive value, there are two real solutions for Eq. (\ref{epsilon_eq}), while it will be shown that only one of them preserves the order of accuracy of the scheme

\begin{equation}
\epsilon_n= \frac{-B_n^* + \sqrt{\hat{\Delta}_n}}{2A_n^*} \quad \text{for} \quad A_n^* \neq 0 \; || \; \hat{\Delta}_n \geq 0.
\end{equation}
On the contrary, when $\hat{\Delta}_n < 0$, there are no real-valued solution for Eq. (\ref{epsilon_eq}). So, $\epsilon_n$ at each step will take one of the following values

\begin{equation}\label{epsilon_sol}
\epsilon_n= \left\{ \begin{array}{lcl}
0 & \mbox{for} & A_n^*=0 \\
\frac{-B_n^* + \sqrt{\hat{\Delta}_n}}{2A_n^*} & \mbox{for} & A_n^* \neq 0 \; || \; \hat{\Delta}_n \geq 0 \\
\mbox{No real value} & \mbox{for} & A_n^* \neq 0 \; || \; \hat{\Delta}_n < 0 \end{array} \right.
\end{equation}

Note that similar to standard RK schemes that conserve linear invariants of the ODE system \cite{hairerGeometricNumericalIntegration2006}, RF-RK schemes also hold this property automatically. 

In the section \ref{Subs_Acc}, it will be proven that, for sufficiently small step sizes, $\hat{\Delta}_n \geq 0$, and a real-valued solution for \ref{epsilon_eq} exists. Moreover, we will show with $\epsilon_n$ set by Eq. (\ref{epsilon_sol}), this new energy conservative scheme has at least the same order of accuracy as the original RK method.

\section{Accuracy and stability of RFRK} \label{Sec_Acc_Stab}

This section demonstrates accuracy preservation of RF-RK schemes in Section \ref{Subs_Acc}, and stability properties of RF-RK schemes compared to unmodified RK schemes in Section \ref{Subs_Stab}.

\subsection{Accuracy} \label{Subs_Acc}

In this section, we show that for sufficiently small step sizes, there exists a real-valued $\epsilon_n$ at each step as a solution for Eq. (\ref{epsilon_order}). Then, it will be proven that starting with a RK scheme of order $p$, the order of accuracy of the RF-RK method will be equal to, or higher than, $p$.

\begin{lemma}
\label{Delta_positive}
For small enough time steps, $\hat{\Delta}_n \geq 0$. So there would be at least one real-valued solution for Eq. (\ref{epsilon_eq})
\end{lemma}

\begin{proof}[Proof of Lemma \ref{Delta_positive}]
As indicated in Section \ref{Subs_RFRK}, when stage derivatives are equal to each other, we have a real-valued solution: $\epsilon_n=0$. However, when stage derivatives are not identical, the existence of real-valued solutions for Eq. (\ref{epsilon_eq}) depends upon the sign of $\hat{\Delta}_n$. Since $\hat{\Delta}_n$ is a function of $A_n^*$, $B_n^*$, and $C_n^*$ (see Eq. (\ref{Delta})), we start with the magnitude of these parameters. 
According to Ketheson \cite{ketchesonRelaxationRungeKutta2019} (proof of Lemma 4), we know that $C^*$ is of order $p-1$

\begin{equation}\label{C_order}
C^*= \mathcal{O}(\Delta t^{p-1}).
\end{equation}

Regarding $B^*$,we first write the Taylor series expansion for stage derivatives $f_j$ up to their linear term \cite{hairerSolvingOrdinaryDifferential1996}

\begin{equation}
f_j= f(t_n+c_j\Delta t, y_j)= f_0 + f_j^{(1)}\Delta t + \mathcal{O}(\Delta t^{2}),
\end{equation}
where $f_0= f(t_n, u^n) $, and $f_j^{(1)}$ represents the linear coefficient within the Taylor series expansion for $f_j$. So, the Taylor series for the inner product of stage derivatives $f_i$ and $f_j$ becomes

\begin{equation}\label{fifj_Taylor}
\langle f_i, f_j \rangle= \langle f_0, f_0 \rangle + \bigg( \langle f_i^{(1)}, f_0 \rangle + \langle f_0, f_j^{(1)} \rangle \bigg) \Delta t + \mathcal{O}(\Delta t^{2}).
\end{equation}
Putting this expression in the definition of $B^*$, Eq. (\ref{B}), the Taylor series expansion for this parameter can be written as
\begin{equation}
B^* = \bigg(  -2\sum_{i,j=1}^s k_ia_{ij} +2 \sum_{i,j=1}^s k_ib_j   \bigg) \langle f_0, f_0 \rangle + \mathcal{O}(\Delta t) .
\end{equation}
The constant portion of this series is composed of two terms. The first, $-2\sum_{i,j} k_ia_{ij}= -2\sum_{i} k_ic_i$, is nonzero by the definition of $k_i$ terms, Eq. (\ref{k_conditions}), while the second is zero, again by the conditions set on $k_i$, $2\sum_{ij} k_ib_j= 2(\sum_{i} k_i= 0)(\sum_{j} b_j=1)=0$. So, we are left with the following Taylor expansion for $B^*$

\begin{equation}
B^* = \bigg( -2 \sum_{i,j=1}^s k_ic_i  \bigg) \langle f_0, f_0 \rangle + \mathcal{O}(\Delta t).
\end{equation}
This means that $B^*$ has a non-zero constant term in its Taylor series

\begin{equation}\label{B_order}
B^*= \mathcal{O}(1).
\end{equation}
Lastly, it can be shown that $A^*$ is of order of a positive integer $m$

\begin{equation}\label{A_order}
A^*= \mathcal{O}(\Delta t^m).
\end{equation}
Therefore, by using the order conditions for $A^*$, $B^*$, and $C^*$, we can write the following for $\hat{\Delta}_n$

\begin{equation}
\begin{split}
\hat{\Delta}_n &= (B_n^*)^2 - 4 A_n^* C_n^*\\
&= (\mathcal{O}(1))^2 - 4 (\mathcal{O}(\Delta t^m))(\mathcal{O}(\Delta t^{p-1})).
\end{split}
\end{equation}
This shows that $\hat{\Delta}_n$ has a positive constant term in its Taylor series, and it will take a positive value if the step size is small enough. Therefore, according to Eq. (\ref{epsilon_sol}), there is at least one real-valued solution for Eq. (\ref{epsilon_eq}) when the step size is sufficiently small.

\end{proof}

Having a positive $\hat{\Delta}_n$ for small step sizes, we can provide a statement about the magnitude of $\epsilon_n$ defined in Eq. (\ref{epsilon_sol}).

\begin{lemma}
\label{epsilon_order}
Having a parent RK method of order p, $\epsilon_n$ in Eq. (\ref{epsilon_sol}) converges to zero with a rate of $p-1$ provided that time steps are sufficiently small
\[  \epsilon_n= \mathcal{O}(\Delta t^{p-1})  . \]  

\end{lemma}

\begin{proof}[Proof of Lemma \ref{epsilon_order}]

When we have $A_n^*=0$, $\epsilon_n$ is identically zero. Moving forward with $A_n^* \neq 0$, we saw that having small enough step sizes, we will have $\hat{\Delta}_n \geq 0$. Since $B^*$ has a non-zero constant term in its Taylor series, for the square root of $\hat{\Delta}_n$ we can write

\begin{equation}
\begin{split}
\sqrt{\hat{\Delta}_n} &= \sqrt{(B_n^*)^2- 4A_n^*C_n^*} \\
& = \sqrt{(B_n^*)^2 - \mathcal{O}(\Delta t^{m+p-1}) } \\
&= \sqrt{\bigg( B_n^* - \mathcal{O}(\Delta t^{m+p-1}) \bigg)^2} \\
&= B_n^* - \mathcal{O}(\Delta t^{m+p-1}).
\end{split}
\end{equation}
Putting this expression in Eq. (\ref{epsilon_sol}), we obtain the magnitude of $\epsilon_n$ as a function of the time step size

\begin{equation}
\begin{split}
\epsilon_n &= \frac{-B_n^* + \sqrt{\hat{\Delta}_n}}{2A_n^*},\\
&= \frac{-B_n^* + \bigg( B_n^* -  \mathcal{O}(\Delta t^{m+p-1}) \bigg)}{\mathcal{O}(\Delta t^{m})}, \\
&= \mathcal{O}(\Delta t^{p-1}),
\end{split}
\end{equation}
which shows $\epsilon_n$ converges to zero with a rate of $p-1$.

\end{proof}

Now that we have obtained the convergence rate of $\epsilon_n$, we can show the order of accuracy of RF-RK schemes.

\begin{theorem}
\label{RFRK_order}
If the parent/original RK method is of order $p$, the RF-RK scheme is of order $p$ or higher.
\end{theorem}

\begin{proof}[Proof of Theorem \ref{RFRK_order}]

Herein, ($n+1$)th step solution out of the parent RK method is denoted as $u^{n+1}$, while the corresponding solution with the RF-RK scheme is denoted $u_{RF}^{n+1}$. Having a parent RK scheme of order $p$ means that the Taylor series for the exact solution $u(t_n+ \Delta t)$ and for $u^{n+1}$ coincide up to the term $\Delta t^p$ \cite{hairerSolvingOrdinaryDifferential1996}

\begin{equation}\label{RK_order}
|| u(t_n+ \Delta t) - u^{n+1} || \leq K \Delta t^{p+1}.
\end{equation}

For $u(t_n+ \Delta t)$, it is possible to write a Taylor series \cite{hairerSolvingOrdinaryDifferential1996} 

\begin{equation}
\begin{split}
u(t_n+ \Delta t) = u(t_n) &+ \Delta tf(t_n, u^n)+ \frac{\Delta t^2}{2} \bigg( f_t+ f_uf   \bigg)(t_n, u^n)  \\
 &+ \frac{\Delta t^3}{6} \bigg( f_{tt}+ 2f_{tu}f + f_{uu}ff + f_uf_t+ f_uf_u f \bigg) (t_n, u^n) + \mathcal{O}(\Delta t^{4}).
\end{split}
\end{equation}
Also, it is possible to write a Taylor series for $u^{n+1}$ 
\begin{equation}\label{u_RK_Taylor}
\begin{split}
u^{n+1}= u(t_n) &+ \Delta t \sum_{i=1}^s b_i f(t_n, u^n) + \frac{ \Delta t^2 }{2} \sum_{i=1}^s 2b_ic_i \bigg( f_t+ f_uf   \bigg)(t_n, u^n) \\
& + \frac{\Delta t^3}{6} \bigg( 3\sum_{i=1}^s b_ic_i^2 ( f_{tt} + 2f_{tu}f +   f_{uu}ff)+ 6\sum_{i,j=1}^s b_ia_{ij}c_j ( f_uf_t + f_uf_uf)  \bigg)(t_n, u^n)  + \mathcal{O}(\Delta t^{4}).
\end{split}
\end{equation}
From this, the following Taylor series for $u_{RF}^{n+1}$ can be obtained easily by replacing $b_i$ with $\hat{b}_i$

\begin{equation}\label{u_RFRK_Taylor}
\begin{split}
u_{RF}^{n+1}= u(t_n) &+ \Delta t \sum_{i=1}^s \hat{b}_i f(t_n, u^n) + \frac{ \Delta t^2 }{2} \sum_{i=1}^s 2\hat{b}_ic_i \bigg( f_t+ f_uf   \bigg)(t_n, u^n) \\
& + \frac{\Delta t^3}{6} \bigg( 3\sum_{i=1}^s \hat{b}_ic_i^2 ( f_{tt} + 2f_{tu}f +   f_{uu}ff)+ 6\sum_{i,j=1}^s \hat{b}_ia_{ij}c_j ( f_uf_t + f_uf_uf)  \bigg)(t_n, u^n)  + \mathcal{O}(\Delta t^{4}).
\end{split}
\end{equation}
Since we know from Eq. (\ref{hat_b}) that $\hat{b}_i - b_i= k_i\epsilon_n$, the difference between the two Taylor series for $u^{n+1}$ and $u_{RF}^{n+1}$ becomes

\begin{equation}\label{u_Taylor_difference}
\begin{split}
u_{RF}^{n+1} - u^{n+1}=&  \Delta t \epsilon_n \bigg[ \sum_{i=1}^sk_i  \bigg] f(t_n, u^n) + \frac{\Delta t^2}{2} \epsilon_n \bigg[ 2 \sum_{i=1}^s k_ic_i \bigg] \bigg( f_t+ f_uf   \bigg)(t_n, u^n) \\
& + \frac{\Delta t^3}{6} \epsilon_n \bigg( 3\sum_{i=1}^s k_ic_i^2 ( f_{tt} + 2f_{tu}f +   f_{uu}ff)+ 6\sum_{i,j=1}^s k_ia_{ij}c_j ( f_uf_t + f_uf_uf)  \bigg)(t_n, u^n)  + \epsilon_n \mathcal{O}(\Delta t^{4}).
\end{split}
\end{equation}
From Eq. (\ref{k_conditions}) we already know that $\sum_{i} k_i=0$ and $\sum_{i} k_ic_i \neq 0$. So, on the right hand side $\epsilon_n= \mathcal{O}(\Delta t^{p-1})$ is multiplied by $\Delta t^m, m \geq 2$. Therefore, the difference between the Taylor series for RF-RK and RK becomes

\begin{equation}
u_{RF}^{n+1} - u^{n+1}= \mathcal{O}(\Delta t^{p+1}),
\end{equation}
and from Eq. (\ref{RK_order} ) we conclude that the order of accuracy of the RF-RK scheme is at least $p$

\begin{equation}
u(t_n + \Delta t) - u_{RF}^{n+1}= \mathcal{O}(\Delta t^{p+1}).
\end{equation}

\end{proof}

\subsection{Stability} \label{Subs_Stab}

In this part, the linear stability of RF-RK schemes with respect to the original RK methods will be examined. Assume the problem of interest satisfies 

\begin{equation}
u'(t)= \lambda u(t),
\end{equation}
$\lambda$ being a complex number. After a step of length $\Delta t$, the exact solution will be multiplied by $e^{z}, z= \lambda \Delta t$. However, the approximate solution by the RK integration scheme will be multiplied by its so-called stability polynomial, $R(z)$ \cite{hairerSolvingOrdinaryDifferential1996}. This stability polynomial also defines the region of linear stability, which is the area where the magnitude of $R(z)$ is less than or equal one. So, we can compare the stability regions from the original RK method and its RF-RK counterpart, by comparing their stability polynomials.

For an RK integration scheme, the stability polynomial can be obtained from \cite{butcherNumericalMethodsOrdinary2008}
\begin{equation}
R(z)= 1 + z b^{T}(I- zA)^{-1} e,
\end{equation}
where $e$ is a vector of ones. In an expanded polynomial form, this will be
\[ R(z)= 1 + z  \sum_{i=1}^s b_i  + z^2  \sum_{i=1}^s b_ic_i   + z^3  \sum_{i,j=1}^s b_ia_{ij}c_j  + ...+ z^s  \sum_{i,j=1}^s b_i g_s,  \]
Where $g_s$ is a function of entries $A$ and $c$ only. By substituting $b_i$ with $\hat{b}_i= b_i + k_i\epsilon_n$ we can create the stability polynomial for an RF-RK scheme at the step $n$, $R_{RF,n}(z)$, as a function of $R(z)$
\[ R_{RF,n}(z) = R(z) + \epsilon_n \bigg( z\sum_{i=1}^s k_i + z^2 \sum_{i=1}^s k_ic_i + z^3 \sum_{i,j=1}^s k_ia_{ij}c_j +...+  z^s \sum_{i=1}^s k_i g_s(A,c) \bigg).  \]
Since we have set $\sum_{i=1}^s k_i=0$, $\epsilon_n$ which is $\mathcal{O}(h^{p-1})$ gets multiplied by $(z=\lambda \Delta t)^m, m \geq 2$. So, we can say the stability polynomials for a RK method and its energy conservative RF-RK scheme at each point $z$ are close to each other with a difference of order $\mathcal{O}(\Delta t^{p+1})$
 \[ |R_{RF,n}(z) - R(z)| =  \mathcal{O}(\Delta t^{p+1}).\]
 
\begin{figure}[h!]
\centering

\begin{subfigure}[b]{0.3\textwidth}
\centering
\includegraphics[width=\textwidth]{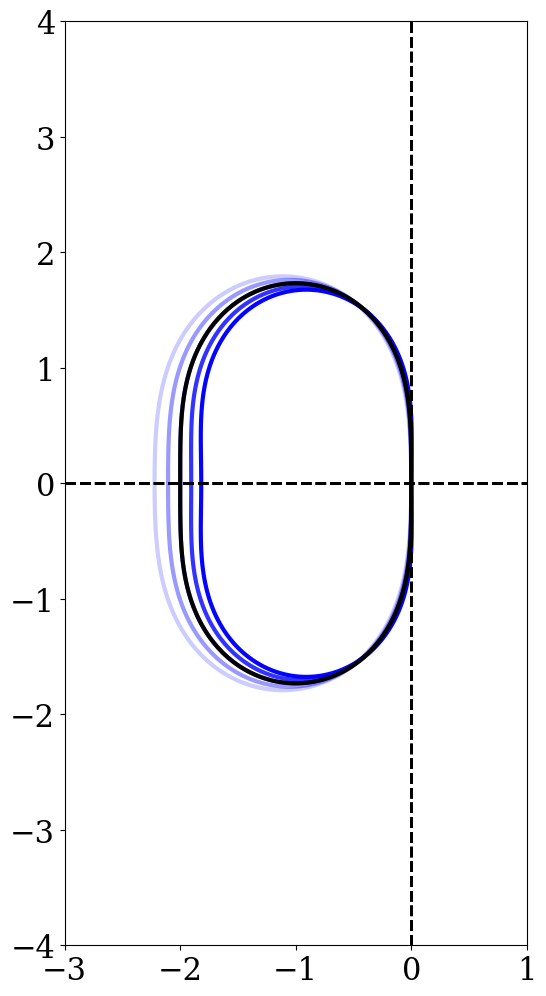}
\caption{2nd order RF-RK method}
\end{subfigure}
\begin{subfigure}[b]{0.3\textwidth}
\centering
\includegraphics[width=\textwidth]{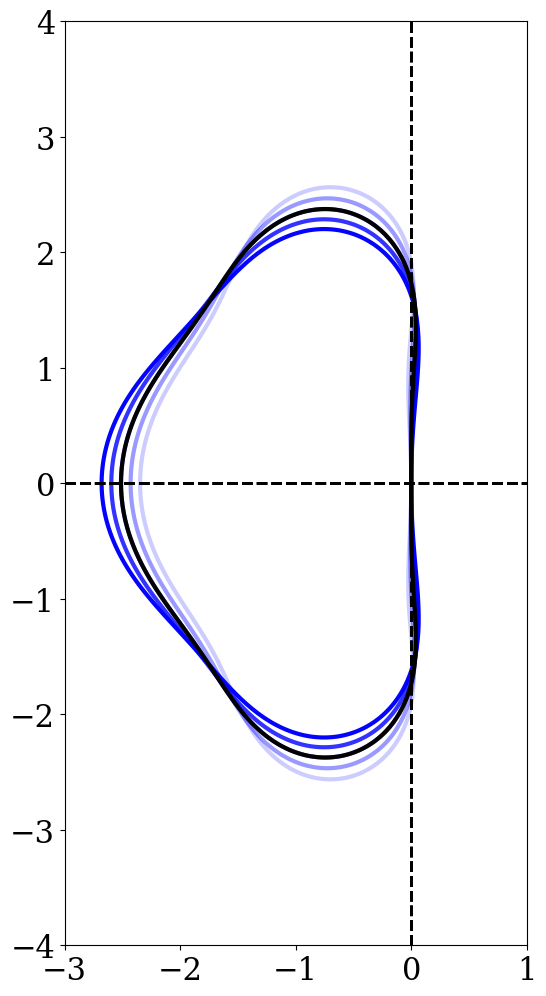}
\caption{ 3rd-order RF-RK method}
\end{subfigure}
\begin{subfigure}[b]{0.3\textwidth}
\centering
\includegraphics[width=\textwidth]{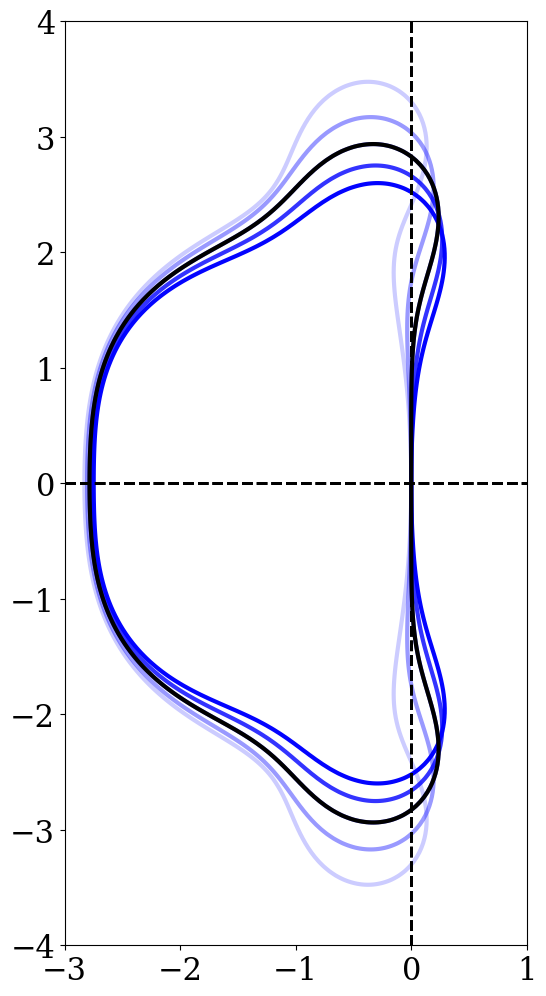}
\caption{ 4th-order RF-RK method}
\end{subfigure}

\caption{Change in linear stability region for integration schemes after application of RF-RK method for $\epsilon_n$ between $-0.05$ to $0.05$, where the black lines belong to $\epsilon_n=0$, and the darker blue lines belong to positive values for $\epsilon_n$. Note that these values for $\epsilon_n$ are chosen large to exhibit its dependence clearly. }
\label{Fig_Stability_region}
\end{figure}

While for non-zero $\Delta t$ the stability limit for the RF-RK method can be different from the original RK method, the two stability limits tend to converge as the step size decreases. To have better insight, Figure \ref{Fig_Stability_region} shows the change in stability regions after applying the RF-RK technique. Note that these plots are obtained with $\epsilon_n$ ranging from $-0.05$ to $0.05$, which are much larger than practical values for $\epsilon_n$ as will be shown later. Clearly, the stability limit along both axes can be changed by the RF-RK method, and this change gets amplified by higher values for $\epsilon_n$. So, it can be expected that by using smaller time steps, $\epsilon_n$ in turn will be smaller in magnitude, and the linear stability region will be modified less.

\section{Numerical examples} \label{Sec_exs}

This section presents the usage of RF-RK schemes for validation test cases. To better compare the results of RF-RK with the already developed R-RK technique, the example cases of Ketcheson \cite{ketchesonRelaxationRungeKutta2019} are followed here. This entails using several types of base RK schemes for different linear and non-linear illustrative problems. 

Types of base RK schemes that are used here are standard RK schemes, SSP methods of  \cite{ketchesonHighlyEfficientStrong2008} called SSPRK, and fifth order method of \cite{bogackiEfficientRungeKuttaPair1996} called BSRK \cite{ketchesonRelaxationRungeKutta2019}. The number of stages and the order of each scheme are indicated as $(s,p)$. For instance, SSPRK(10,4) is the 4th-order SSPRK method with 10 stages. Using this notation, the base integration schemes used here based on Kecheson \cite{ketchesonRelaxationRungeKutta2019} are

\begin{itemize}
\item SSPRK(2,2) \cite{ketchesonHighlyEfficientStrong2008},
\item SSPRK(3,3) \cite{ketchesonHighlyEfficientStrong2008},
\item RK(4,4)    \cite{sunStabilityFourthOrder2017},
\item BSRK(8,5)  \cite{ketchesonRelaxationRungeKutta2019}.
\end{itemize}
We will use the term "R" or "RF" before the names of base methods to represent their relaxation or relaxation-free energy conservative ones, respectively. For example, RF-RK(4,4) is the relaxation-free version of RK(4,4). 

In the following examples, the $k$ multipliers which are needed to employ RF-RK technique are chosen to be $[1,-1]$ for RF-SSPRK(2,2), $[2,-1,-1]$ for RF-SSPRK(3,3), $[1,2,-2,-1]$ for RF-RK(4,4), and $[2,-1,-1,0,0,0,0,0]$ for RF-BSRK(8,5), although other consistent values may also be used.

\subsection{A normal, linear, autonomous problem}

The first example case from \cite{ketchesonRelaxationRungeKutta2019} concerns the application of the Fourier spectral collocation technique to discretize space in a 1D linear advection problem. The domain, which is periodic with a length of $2\pi$, is discretized with $m=128$ points. The resulting semi-discrete ODE system becomes
  
\begin{equation}\label{ex1_ode}
u'(t)= -Du(t),
\end{equation}
where $D$ is the $m \times m$ skew-Hermitian Fourier spectral differentiation matrix. Since this matrix is normal, we can obtain the biggest linearly stable step size using its eigenvalues, called $\lambda$. For this problem, all the eigenvalues lay on the imaginary axes
\[ \lambda= \pm ik, \; k= \frac{m}{2}-1, \frac{m}{2}-2, .., 0   . \]   

To have linear stability, all of these eigenvalues should fall within the stability region of the employed integration scheme. By using a base RK scheme with an imaginary axis stability limit of $I(A,b)$, the biggest stable step size, $\Delta t_{max}$, becomes

\[ \Delta t_{max}= \frac{I(A,b)}{max|\lambda|}=  \frac{I(A,b)}{m/2-1}.  \]
Note that relative to reference \cite{ketchesonRelaxationRungeKutta2019}, a minus one is added to the denominator as the modes are mirrored about the $\lambda=0$ mode. \par
From the section \ref{Subs_Stab}, we saw that by modifying the $b$ parameters to $\hat{b}$ through the RF-RK technique, the stability limit would change, which in turn affects the biggest linearly stable step size. This change will also happen for the R-RK method \cite{ketchesonRelaxationRungeKutta2019}. Nevertheless, in this exercise the step size limit for the base RK scheme will be used as a reference, and each step size is defined using a multiplier $\mu$, where $\Delta t= \mu \Delta t_{max}$.

To study the energy change over time, one may decompose the solution vector into its spectral components and see the change in amplitude of each mode over time. Spectral decomposition can be performed by discrete Fourier transformation, which gives the following for the $n$th step solution

 \[ u_j^n= \sum_{k=0}^{m-1}\hat{u}_k^n e^{ikx_j} ,  \]
where $\hat{u}_k^n$ is a complex multiplier whose magnitude is the the amplitude of mode $k$ at the time step $n$. With the exact solution, the amplitude of each mode remains constant over time, so the total energy remains unchanged. However, numerical integration may result in dissipation or amplification of each mode, which can lead to a change in total energy. The change in amplitude of each mode $k$ after $n$ time steps can be represented by their relative amplification factor \cite{ketchesonRelaxationRungeKutta2019} 
\[ \frac{|\hat{u}_k^n|- |\hat{u}_k^0|}{|\hat{u}_k^0|}.   \]
Note that since $\hat{u}_k^n$ and $\hat{u}_{m-k}^n$ are complex conjugate and their magnitudes are equal, it is sufficient to study the amplitude change for half of the modes.

\begin{figure}[!htb]
\centering
\includegraphics[width=0.8\textwidth]{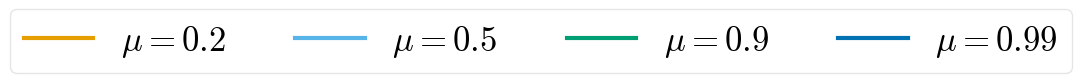}

\begin{subfigure}[b]{0.3\textwidth}
\centering
\includegraphics[width=\textwidth]{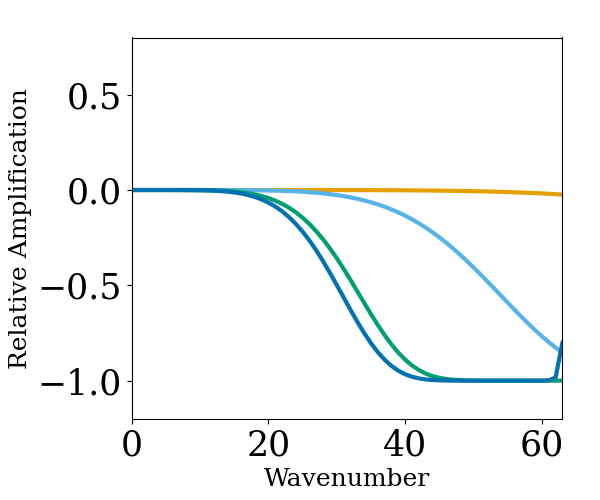}
\caption{RK(4,4) method}
\end{subfigure}
\begin{subfigure}[b]{0.3\textwidth}
\centering
\includegraphics[width=\textwidth]{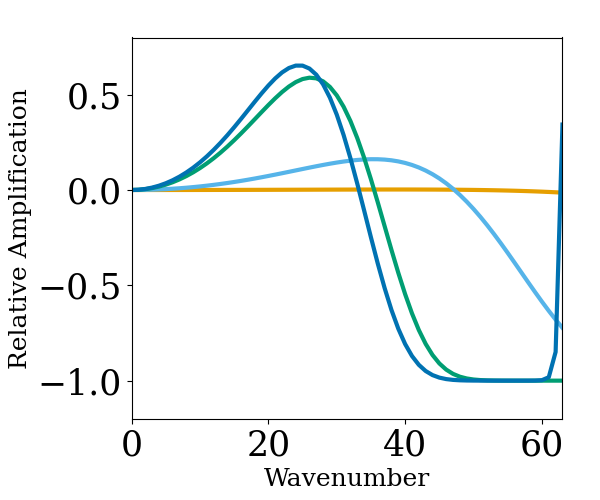}
\caption{R-RK(4,4) method}
\end{subfigure}
\begin{subfigure}[b]{0.3\textwidth}
\centering
\includegraphics[width=\textwidth]{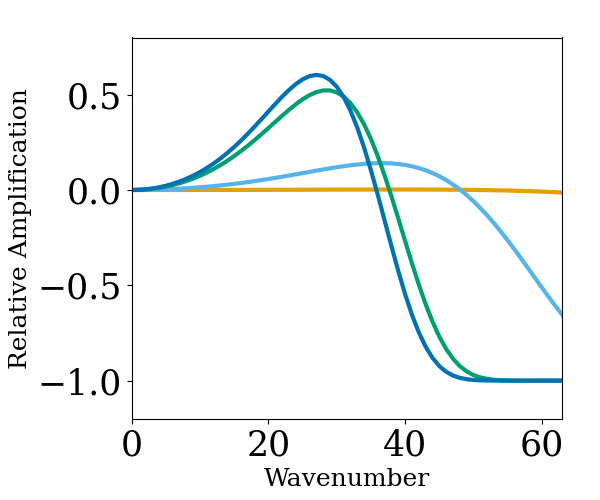}
\caption{RF-RK(4,4) method}
\end{subfigure}

\caption{Relative amplification for each wavelength when using white noise initial data for the problem \ref{ex1_ode}. It is integrated up to $ t_f=1$, with the step sizes $\Delta t= \mu \Delta t_{max}$. While with the RK(4,4) the total energy has decreased, the two energy conserving schemes amplifyed lower modes to cancel out the energy loss by high wavenumber modes.}
\label{Fig_ex1_wave_noise}
\end{figure}

For the first case, the problem is initialized with a white noise input ($\hat{u}_k^0= e^{i\theta_k}$, $\theta_k$ being random) to distribute energy equally across all wave numbers. Then, we record the amplitude change for each mode after a final time of $t_f=1$, for each integration scheme and time step size. The results are depicted in Figure \ref{Fig_ex1_wave_noise}, where the output for standard RK(4,4) is compared with its two energy conserving counterparts. With the RK(4,4) scheme, high wave number modes exhibit significant damping when increasing the step size. So, the total energy will decrease more by using larger step sizes. On the other hand, with both energy conserving schemes, energy loss at high wave numbers is compensated for amplification of lower modes, such that the total energy remains unchanged for each step size. Moreover, a slight difference in stability limits for energy conserving schemes and the original method is visible in the behavior of high wave number modes with the time step $\Delta t= 0.99 \Delta t_{max}$. This time step is slightly larger than the stability limit for R-RK method, which is why the highest modes started to amplify sharply, overpassing $0$. However, this step size is smaller than the stability limit for RF-RK method, and there is no such sharp amplification. Finally, we note the similarity of the R-RK and RF-RK results, while the RF-RK scheme did not require relaxation of the step size.

\begin{figure}[!htb]
\centering
\includegraphics[width=0.8\textwidth]{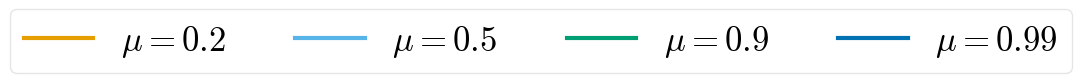}

\begin{subfigure}[b]{0.3\textwidth}
\centering
\includegraphics[width=\textwidth]{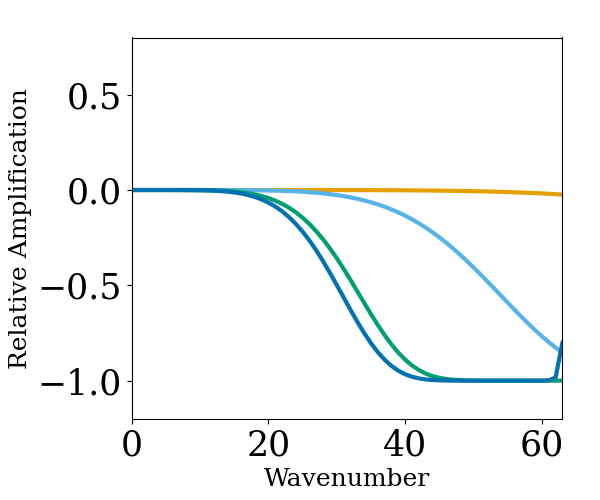}
\caption{RK(4,4) method}
\end{subfigure}
\begin{subfigure}[b]{0.3\textwidth}
\centering
\includegraphics[width=\textwidth]{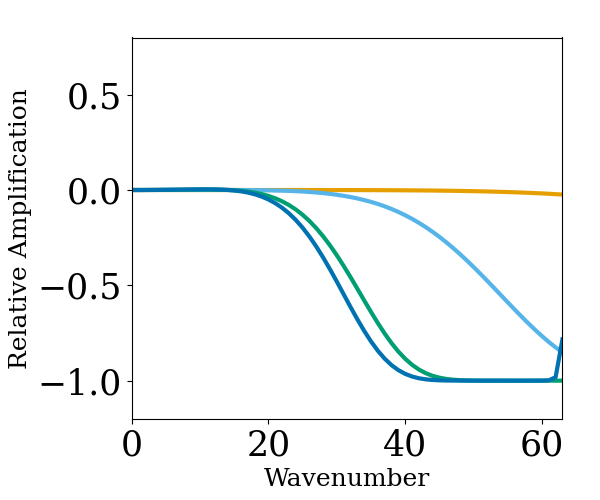}
\caption{R-RK(4,4) method}
\end{subfigure}
\begin{subfigure}[b]{0.3\textwidth}
\centering
\includegraphics[width=\textwidth]{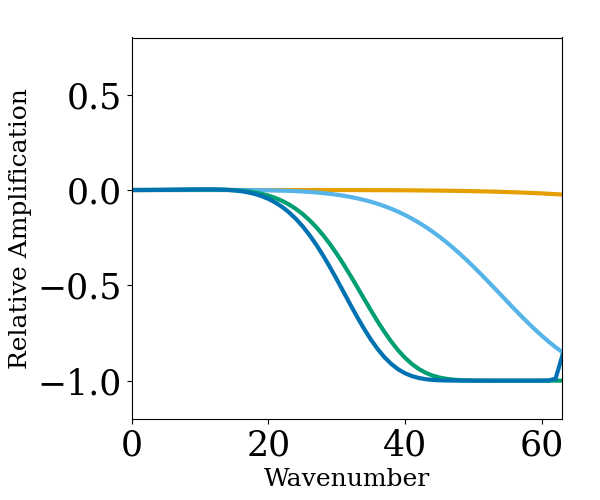}
\caption{RF-RK(4,4) method}
\end{subfigure}

\caption{Relative amplification for each wavelength with smooth initial data for the problem \ref{ex1_ode}. Integration is performed up to $t_f=1$, with the step sizes $\Delta t= \mu \Delta t_{max}$. Since most of the initial energy is distributed among the low wavenumber modes, energy change with RK(4,4) is small, and the curves for RK(4,4) and the energy conserving schemes are similar. }
\label{Fig_ex1_wave_smooth}
\end{figure}

For the second case in this example, we instead distribute the initial energy mainly among the low wavenumber modes by employing a smooth initial data

\begin{equation}
U^0= sech^2(7.5(x+1)).
\end{equation}
Figure \ref{Fig_ex1_wave_smooth} shows the amplitude change of each mode with this smooth initial condition. For RK(4,4), low wavenumber modes, which contain most of the  energy, see a negligible change in their amplitude. Therefore, the change in total energy is very small, and the conservative schemes need to impose little modification to cancel out energy loss. This results in similar behavior for energy conserving and non-conserving schemes.

\begin{figure}[!htb] 
\centering
\begin{subfigure}[b]{0.4\textwidth}
\centering
\includegraphics[width=\textwidth]{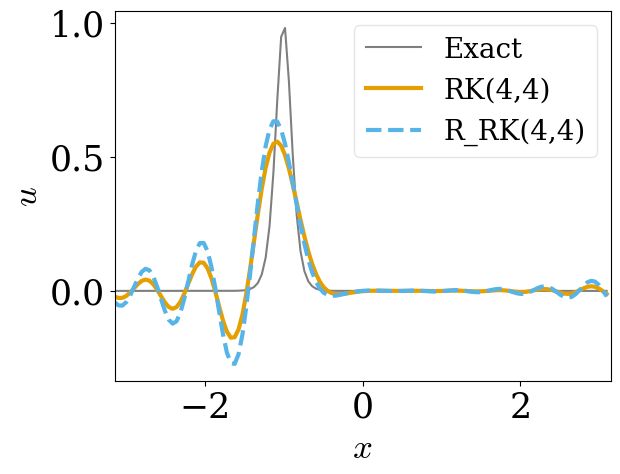}
\caption{RK(4,4) vs. R-RK(4,4) }
\end{subfigure}
\begin{subfigure}[b]{0.4\textwidth}
\centering
\includegraphics[width=\textwidth]{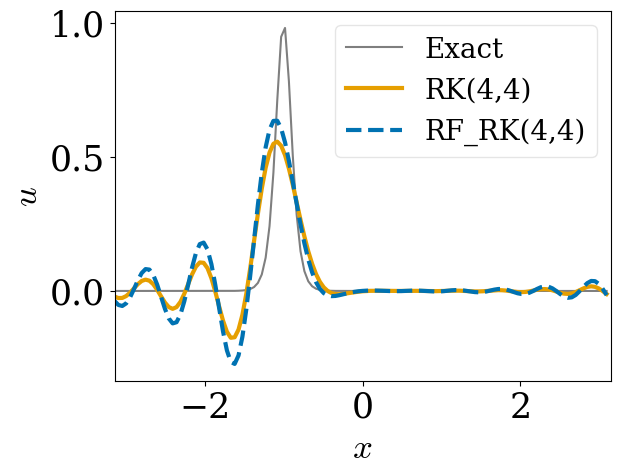}

\caption{ RK(4,4) vs. RF-RK(4,4)}
\end{subfigure}

\caption{Final solutions for problem \ref{ex1_ode} with the smooth initial data, integrated up to $t_f=400\pi$, with $\mu=0.99$. The two energy conserving schemes, R-RK(4,4) and RF-RK(4,4) behave similarly for this step size.}
\label{Fig_ex1_smooth1}
\end{figure}

\begin{figure}[!htb]
\centering
\begin{subfigure}[b]{0.4\textwidth}
\centering
\includegraphics[width=\textwidth]{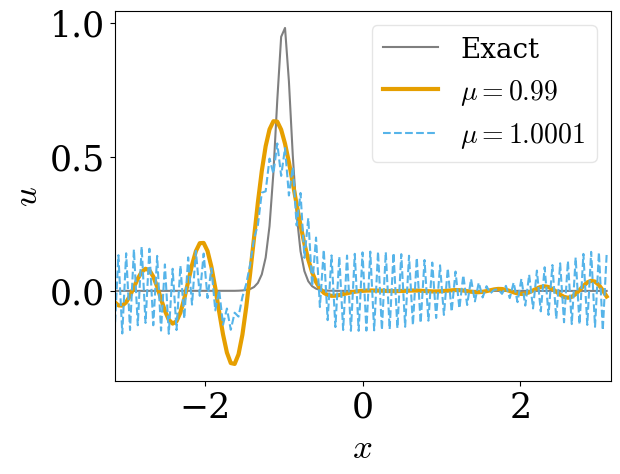}
\caption{ R-RK(4,4) }
\end{subfigure}
\begin{subfigure}[b]{0.4\textwidth}
\centering
\includegraphics[width=\textwidth]{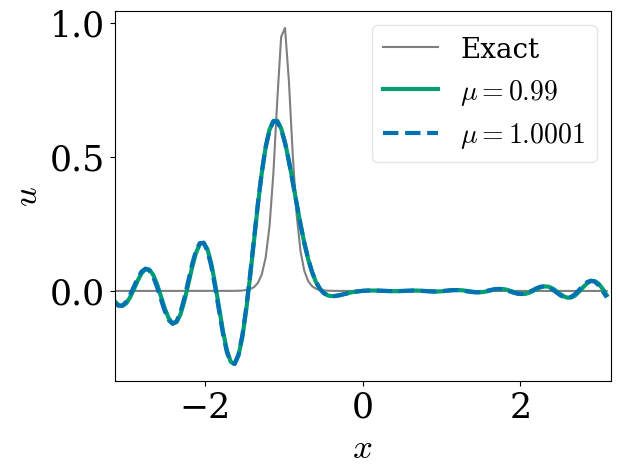}
\caption{ RF-RK(4,4)}
\end{subfigure}

\caption{Comparison of the solutions out of the two energy conserving schemes for the problem \ref{ex1_ode} when it is integrated up to $t_f=400\pi$ with $\mu= 0.99$, and $\mu= 1.0001$. While $\Delta t= 1.0001 \Delta t_{max}$ is slightly higher than the stability limit for R-RK(4,4), it is still within the stability limit for RF-RK(4,4).}

\label{Fig_ex1_smooth2}
\end{figure}

Suggested by \cite{ketchesonRelaxationRungeKutta2019}, the solution for the smooth data after a final time of $t_f= 400 \pi$ with $\mu= 0.99$ is provided in Figure \ref{Fig_ex1_smooth1}. As indicated, for smooth initial data the energy conservative schemes behave similar to the unmodified RK(4,4) scheme. Moreover, for this time step size, the solutions from R-RK(4,4) and RF-RK(4,4) are nearly indistinguishable. Note that $\epsilon_n$ for RF-RK(4,4) in this case remains less than \num{1.25e-3}. For a slightly larger step size when $\mu= 1.0001$, Figure \ref{Fig_ex1_smooth2} compares the behavior of the two energy conserving schemes. While R-RK(4,4) sees linear instability with high amplitude oscillations, RF-RK(4,4) still remains stable. Again for this step size, $\epsilon_n$ for RF-RK(4,4) remains less than \num{1.25e-3}. With larger step sizes, the relaxation-free scheme may not find a real solution for $\epsilon_n$, and the relaxation technique encounters problems in completing the simulation because $\gamma_n$ tends to zero. However, this is not unexpected, as both have exceeded their linear stability limits.

\subsection{A linear energy-decaying system }

\begin{figure}[ht]
\centering

\includegraphics[width=0.5\textwidth]{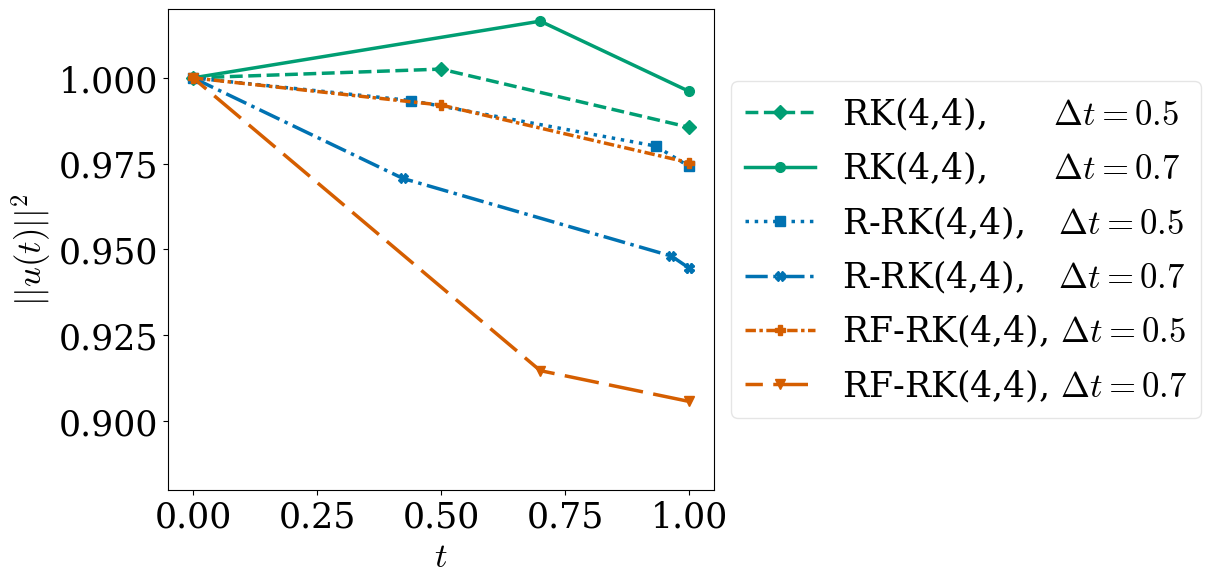}
\caption{ Change in energy after one time step for a dissipative problem integrated with RK(4,4) and the two monotonicity-preserving techniques, using step sizes  of $\Delta t=0.5$ and $\Delta t=0.7$. With RK(4,4) energy increased at the first step, but R-RK and RF-RK methods preserved the monotonicity of the problem. However, R-RK has modified the input step size, while for RF-RK it is constant.}

\label{Fig_ex2}

\end{figure}

Sun \& Shu \cite{sunStabilityFourthOrder2017} demonstrated that for a dissipative semidiscrete system, integrating with RK(4,4) may adversely lead to an increase in energy after the first integration step, no matter how small the step size is. An indicative example is a linear dissipative system in the form of $u'(t)=L u(t)$

\begin{equation}
\begin{bmatrix}
u_1  \\
u_2 \\
u_3
\end{bmatrix}'= \begin{bmatrix}
-1 & -2 & -2 \\
0 & -1 & -2 \\
0 & 0 & -1
\end{bmatrix} \begin{bmatrix}
u_1 \\
u_2 \\
u_3
\end{bmatrix},
\end{equation}
with an initial condition equal to the first right singular vector of $R(0.5L)$, with $R(z)$ being the stability polynomial of RK(4,4). Ketcheson \cite{ketchesonRelaxationRungeKutta2019} showed that while RK(4,4) cannot preserve monotonicity for this problem, application of the relaxation technique allows it to be monotonicity preserving. Figure \ref{Fig_ex2} shows energy change after one step through each scheme, while using two different step sizes: $\Delta t= 0.5$ and $\Delta t= 0.7$. As stated, with the standard RK(4,4) energy increases for both step sizes. In contrast, both R-RK(4,4) and RF-RK(4,4) preserve monotonicity and make the energy decrease for each step size. This monotonicity-preservation however comes with the cost of step size modification for the R-RK technique. Table \ref{Tab: 2} demonstrates these step size modifications. With $\Delta t=0.5$, the actual first step size for R-RK becomes $\gamma_1 \Delta t \simeq 0.44$, and for the larger step $\Delta t=0.7$ it becomes $\gamma_1 \Delta t \simeq 0.42$, even smaller than before. However, using the RF-RK technique we achieve monotonicity-preservation without any step size modification.

\begin{table}[H] 
\centering
\caption{Comparison of assigned step sizes and actual step sizes for each monotonicity-preserving scheme. While for R-RK there is a noticeable difference between the input and actual step size, the RF-RK method keeps the step size unchanged.  }
\begin{tabular}{ p{3cm} p{3cm} p{3cm} p{3cm} }
\hline
 Integration scheme & Assigned $\Delta t$ & Actual $\Delta t$ & Relative $\Delta t$ change  \\
\hline
R-RK(4,4) & $0.50$ & $0.44$ & $12 \%$ \\
\hline
R-RK(4,4) & $0.70$ & $0.42$ & $40 \%$ \\
\hline
RF-RK(4,4) & $0.50$ & $0.50$ & $0 \%$ \\
\hline
RF-RK(4,4) & $0.70$ & $0.70$ & $0 \%$ \\
\hline
\end{tabular}
\label{Tab: 2}
\end{table}

\subsection{A nonlinear oscillator}

Following another example from the work of Ketcheson \cite{ketchesonRelaxationRungeKutta2019}, we test the integration schemes on a conservative, nonlinear oscillator problem

\begin{equation}
\begin{bmatrix}
u_1 \\
u_2
\end{bmatrix}' = \frac{1}{||u||^2}\begin{bmatrix}
-u_2 \\
u_1
\end{bmatrix}, \; \begin{bmatrix}
u_1(0) \\
u_2(0)
\end{bmatrix}= \begin{bmatrix}
1\\
0
\end{bmatrix},
\end{equation}
which has the following exact analytical solution

\[ \begin{bmatrix}
u_1(t) \\
u_2(t)
\end{bmatrix} = \begin{bmatrix}
\text{cos}(t)\\
\text{sin}(t) \\
\end{bmatrix}. \]

\begin{figure}\label{Figure_3_evol}
\centering
\includegraphics[width=\textwidth]{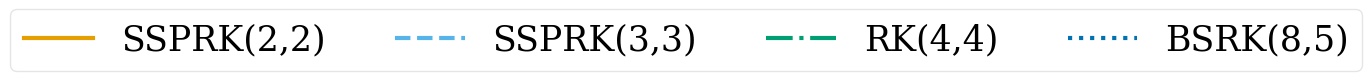}

\begin{subfigure}[b]{0.3\textwidth}
\centering
\includegraphics[width=\textwidth]{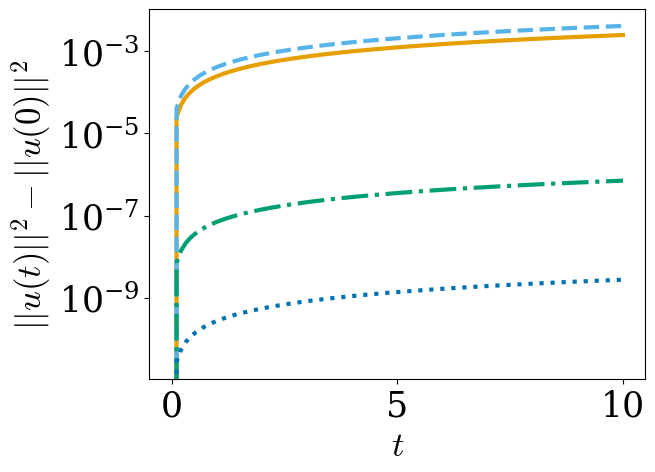}
\caption{ RK }
\label{Figure_3_evol_a}
\end{subfigure}
\begin{subfigure}[b]{0.3\textwidth}
\centering
\includegraphics[width=\textwidth]{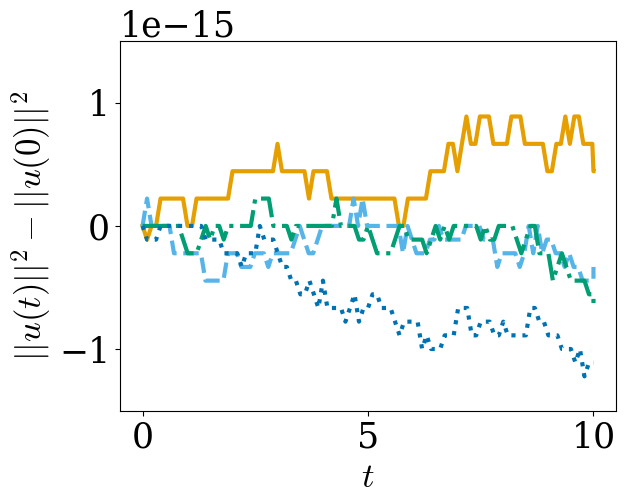}
\caption{R-RK}
\label{Figure_3_evol_b}
\end{subfigure}
\begin{subfigure}[b]{0.3\textwidth}
\centering
\includegraphics[width=\textwidth]{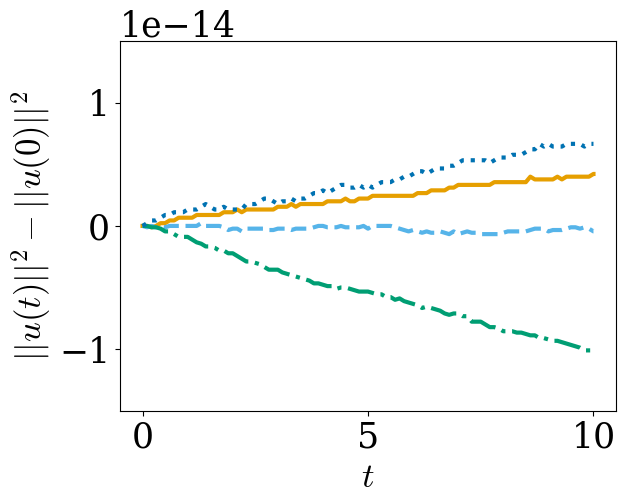}
\caption{RF-RK}
\label{Figure_3_evol_c}
\end{subfigure}

\caption{Evolution of energy for the nonlinear oscillator problem, integrated with different integration schemes and a time step of $\Delta t=0.1$. All the original RK methods increase energy up to their truncation error, but the R-RK and RF-RK techniques conserve energy up to machine precision, and the RF-RK method does so with a constant step size. }
\end{figure}

Figure \ref{Figure_3_evol_a} shows that employing each of the unmodified RK schemes with a time step size of $\Delta t=0.1$ leads to a monotonic increase in energy. On the other hand, presented in Figures \ref{Figure_3_evol_b} and \ref{Figure_3_evol_c}, both R-RK and RF-RK techniques made these RK schemes conserve energy up to machine precision. Note that while using the R-RK method the actual time step ($\gamma_n \Delta t$) is not exactly $0.1$, but it is in the range of $0.0995 \leq \gamma_n \Delta t \leq 0.1 $. On the other hand, with the RF-RK method the step size maintains exactly $0.1$, while the variable $\epsilon_n$ remains in the range of $-0.0015 \leq \epsilon_n \leq 0$.

\begin{figure}
\centering
\includegraphics[width=0.6\textwidth]{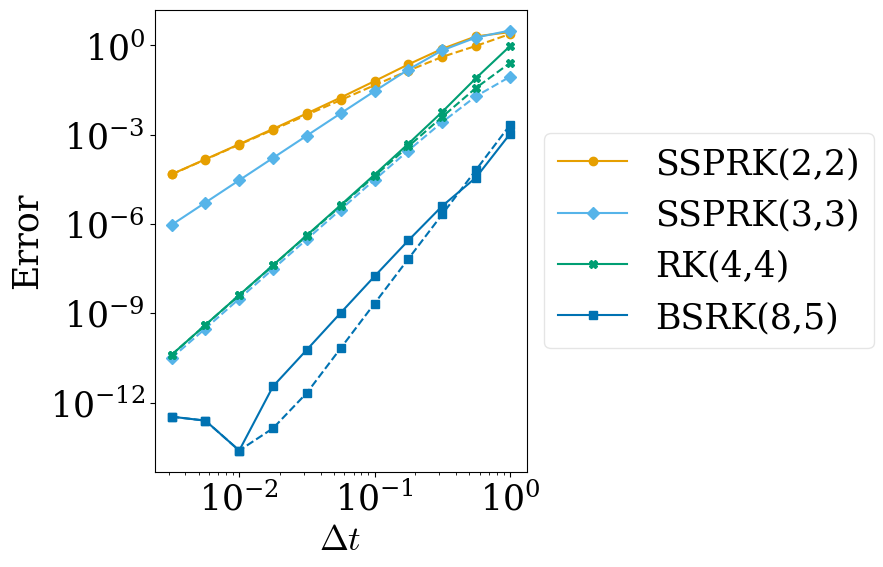}

\caption{Convergence study for the nonlinear oscillator problem. For base RK methods convergence is demonstrated by solid lines, while for their relaxation-free counterpart it is depicted by dashed lines. With RF-RK, the order of accuracy is equal to, or higher than, the corresponding base RK scheme. }
\label{Figure_3_convergence}
\end{figure}

Concerning their accuracy, Figure \ref{Figure_3_convergence} shows solution convergence for unmodified schemes with solid lines and the corresponding RF-RK schemes in dashed lines. It confirms that with the  RF-RK method, the order of accuracy is either the same, or higher, than the base RK method.

\subsection{Burgers' equation}

The last example case which has been used for R-RK methods in \cite{ketchesonRelaxationRungeKutta2019} is an inviscid Burger's problem on a periodic interval of $-1 \leq x \leq 1$ 

\begin{equation}
U_t+ \frac{1}{2}(U^2)_x=0  ,
\end{equation}

\[ U(x,0)= e^{-30x^2}.  \] 
This problem can be transformed to an energy conservative ODE system by discretizing the domain with $50$ equally-spaced points and using the second-order accurate symmetric flux \cite{tadmorEntropyStabilityTheory2003}

\begin{equation}
u_i'(t)= -\frac{1}{\Delta x} (F_{i+1/2} - F_{i-1/2}) , \qquad F_{i+1/2}=\frac{u_i^2 + u_iu_{i+1} + u_{i+1}^2}{6}.
\end{equation}

\begin{figure}[t] 
\centering
\includegraphics[width=\textwidth]{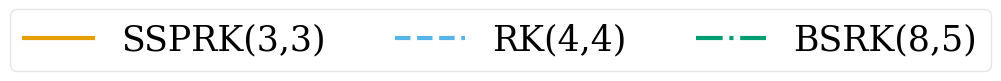}

\begin{subfigure}[b]{0.3\textwidth}
\centering
\includegraphics[width=\textwidth]{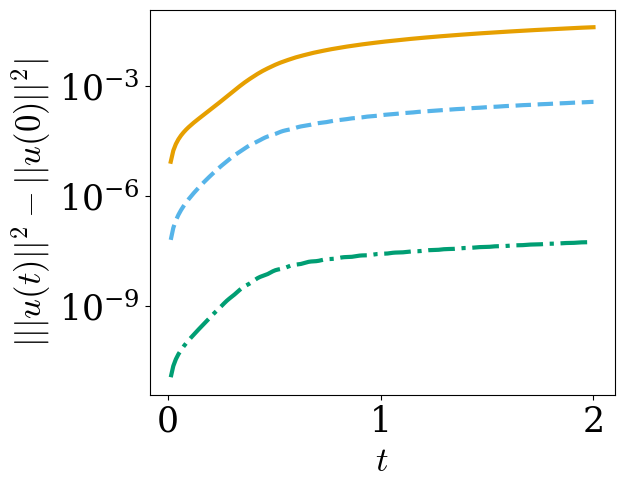}
\caption{ RK }
\end{subfigure}
\begin{subfigure}[b]{0.3\textwidth}
\centering
\includegraphics[width=\textwidth]{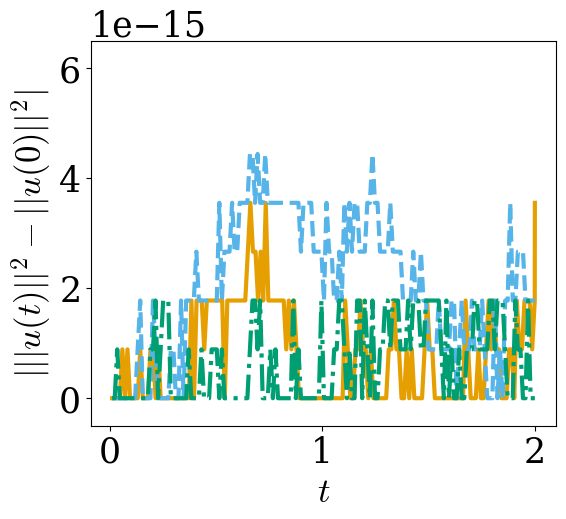}
\caption{R-RK}
\end{subfigure}
\begin{subfigure}[b]{0.3\textwidth}
\centering
\includegraphics[width=\textwidth]{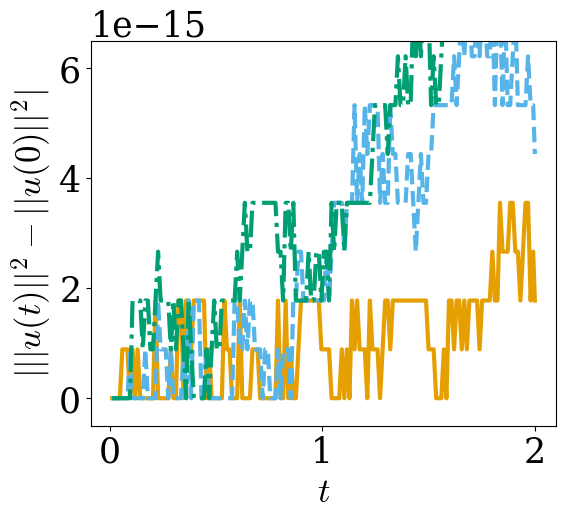}
\caption{RF-RK}
\end{subfigure}

\caption{Evolution of energy for the Burgers' equation with different integration schemes with a step size of $\Delta t=0.3\Delta x$ up to a final time of $t_f=2$. With each of the unmodified RK schemes energy increased monotonically, while energy conserving counterparts conserved energy up to machine precision}
\label{Figure_4_evol}
\end{figure}

\begin{figure}[t] 
\centering

\includegraphics[width=\textwidth]{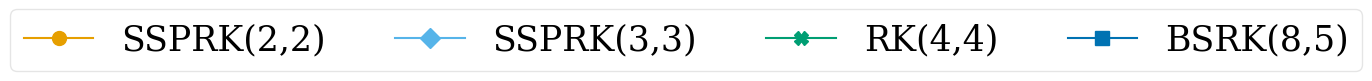}

\begin{subfigure}[b]{0.3\textwidth}
\centering
\includegraphics[width=\textwidth]{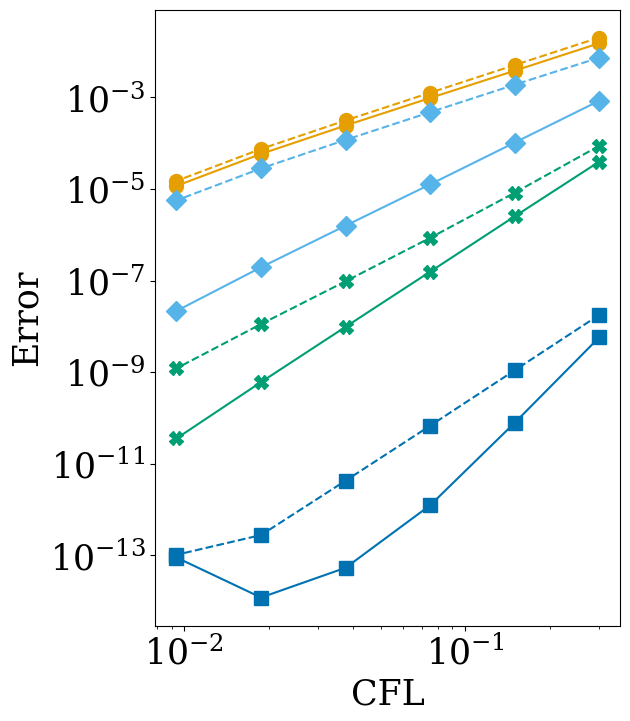}
\caption{ IDT-RK }
\end{subfigure}
\begin{subfigure}[b]{0.3\textwidth}
\centering
\includegraphics[width=\textwidth]{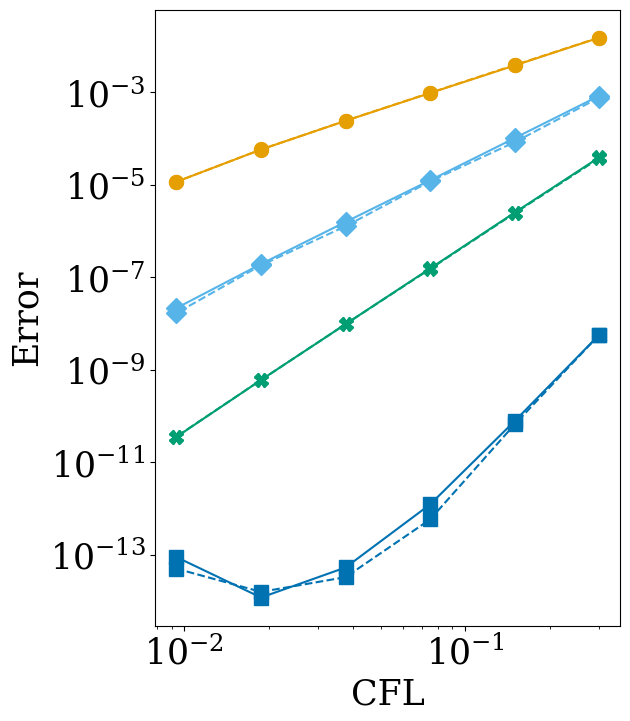}
\caption{R-RK}
\end{subfigure}
\begin{subfigure}[b]{0.3\textwidth}
\centering
\includegraphics[width=\textwidth]{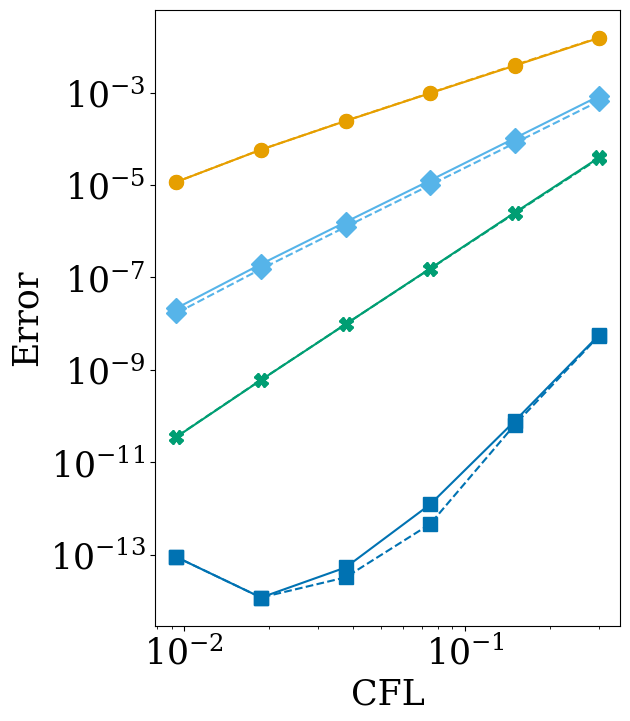}
\caption{RF-RK}
\end{subfigure}

\caption{Convergence analysis for Burger's equation integrated up to a final time of $t_f=0.2$ with the time steps: $\Delta t= \text{CFL} \times \Delta x$. It confirms that while IDT-RK decreases the order of accuracy by one, both R-RK and RF-RK methods retain the order of accuracy of original RK method.}
\label{Fig_4_convergence}
\end{figure}

This problem is integrated with a time step of $\Delta t=0.3\Delta x$ up to a final time of $t_f=2$ using different base RK schemes and their energy conservative counterparts. Figure \ref{Figure_4_evol} shows that with the original RK schemes, there is a noticeable change in the energy of the system, while with both R-RK and RF-RK schemes, energy is conserved up to machine precision. Then, to visualize the order of accuracy of the schemes, convergence analysis has been performed at $t_f=0.2$ with a fixed spatial discretization and different input time steps: $\Delta t= \text{CFL} \times \Delta x$, $\text{CFL}= 0.3 \times 0.5^{0, 1, ..,6}$. The results in Figure \ref{Fig_4_convergence} confirm that with R-RK order of accuracy is retained, but with IDT-RK order of accuracy decreases by one. RF-RK however, benefits from a fixed step size, similar to IDT-RK, and preserved accuracy, similar to R-RK schemes.

\section{Conclusions} \label{Sec_conclusion}

We have proposed a new family of energy conservative Relaxation-Free Runge-Kutta schemes. These schemes introduce a simple modification to the Butcher tableau coefficients that allow conservation of energy, preservation of the order of accuracy of the base RK scheme, while maintaining a constant time-step size. This is in contrast to classical RK schemes, which are not energy conservative, IDT-RK schemes, which do not maintain their base scheme's order of accuracy, and R-RK schemes, which do not maintain a constant step size. In this sense, the proposed RF-RK schemes are superior to these previous schemes in several respects. Numerical results demonstrate that these aforementioned properties are observed in practice for a normal, linear, autonomous problem, a linear energy-decaying system, a non-linear oscillator, and Burgers equation. The RF-RK schemes consistently conserved energy while maintaining a constant step size and the order of accuracy of their base RK scheme.

The proposed RK-RK schemes present a promising framework for non-linear stability of time dependent systems of PDEs. Future work will focus on extension to other forms of entropy, and application to systems with multiple such constraints.

\section*{Acknowledgments}

The authors acknowledge financial support from the Natural Sciences and Engineering Research Council of Canada (NSERC) and the Fonds de Recherche du Québec - Nature et Technologies (FRQNT) via the NOVA program. Additionally, the authors would like to thank Siva Nadarajah, Alexander Cicchino, and Carolyn Pethrick for their helpful discussions at the early stages of this work.

\bibliographystyle{unsrt}

\bibliography{Time_integration}

\end{document}